\documentclass[12pt,fleqn]{iopart}
\usepackage[english]{babel}
\usepackage{amssymb}
\usepackage{amsfonts}
\usepackage{euscript}

\newcommand{\textfrac}[2]{{\textstyle\frac{#1}{#2}}}

\begin{document}
\title[Symmetry Pseudo-Group and Covering of Second Heavenly Equation]{%
Maurer--Cartan Forms of the Symmetry Pseudo-Group and the Covering of Pleba\~nski's 
Second Heavenly Equation
}

\author{Oleg I. Morozov}

\address{Department of Mathematics, Moscow State Technical University of Civil Aviation, Kronshtadtskiy Blvd 20, Moscow 125993, Russia
\\
oim{\symbol{64}}foxcub.org}

\begin{abstract}
We derive Wahlquist--Estabrook forms of the covering of Pleba\~nski's second heavenly equation from Maurer--Cartan forms of its symmetry pseudo-group.
%\keywords{Lie pseudo-groups \and Maurer--Cartan forms \and symmetries of differential equations 
%\and coverings of differential equations}
\end{abstract}

\ams{58H05, 58J70, 35A30}

\maketitle

\section{Introduction}
In our preceding papers \cite{Morozov2007} -- \cite{Morozov2009} it was shown that for a number of nonlinear partial  differential equations ({\sc pde}s) with three independent variables Wahlquist--Es\-ta\-brook forms of their coverings can be derived from Maurer--Cartan forms of their symmetry pseudo-groups. In this paper we consider Pleba\~nski's second heavenly equation, \cite{Plebanski},
\begin{equation}
u_{xz}=u_{ty}+u_{xx}\,u_{yy}-u_{xy}^2,
\label{PlebanskiEquation}
\end{equation}
describing self-dual metrics in theory of gravitation. This equation can be obtained  as the compatibility condition for the following system of {\sc pde}s, \cite{Husain,BogdanovKonopelchenko}, cf. \cite[Eq. (3.13)]{Plebanski}:
\begin{equation}
q_t = (u_{xy} - \lambda) \, q_x - u_{xx}\,q_y,
\qquad
q_z = u_{yy} \, q_x - (u_{xy} + \lambda)\,q_y,
\label{covering}
\end{equation}
where $\lambda$ is an arbitrary constant. This condition is equivalent to the commutativity of the following four infinite-dimensional vector fields
\begin{eqnarray}
\widetilde{D}_t &=& \bar{D}_t 
+ \sum \limits_{i,j \ge 0} \widetilde{D}^i_x \widetilde{D}^j_y \left(
(u_{xy} - \lambda) \, q_{1,0} - u_{xx}\,q_{0,1}
\right)\,\frac{\partial}{\partial q_{i,j}},
\nonumber
\\
\widetilde{D}_x &=& \bar{D}_x 
+ \sum \limits_{i,j \ge 0} q_{i+1,j}\,\frac{\partial}{\partial q_{i,j}},
\\
\widetilde{D}_y &=& \bar{D}_y 
+ \sum \limits_{i,j \ge 0} q_{i,j+1}\,\frac{\partial}{\partial q_{i,j}},
\\
\widetilde{D}_z &=& \bar{D}_z 
+ \sum \limits_{i,j \ge 0} \widetilde{D}^i_x \widetilde{D}^j_y \left(
u_{yy} \, q_{1,0} - (u_{xy} + \lambda)\,q_{0,1}
\right)\,\frac{\partial}{\partial q_{i,j}},
\nonumber
\end{eqnarray}
where $\bar{D}_t$, $\bar{D}_x$, $\bar{D}_y$ and $\bar{D}_z$ are restrictions of the total derivatives $D_t$, $D_x$, $D_y$ and $D_z$ to the infinite prolongation of (\ref{PlebanskiEquation}).  This construction is called a {\it covering}, \cite{KV84} -- \cite{KV99}. Dually coverings can be defined by means of differential 1-forms called {\it Wahlquist--Estabrook forms}, \cite{WE}. For Eq. (\ref{PlebanskiEquation}) an ideal of the Wahlquist--Es\-ta\-brook forms is ge\-ne\-ra\-ted by the following forms:
\begin{eqnarray}
\omega_{0,0} &=& 
dq_{0,0}
-((u_{xy}-\lambda)\,q_{1,0}-u_{xx}\,q_{0,1})\,dt
-q_{1,0}\,dx - q_{0,1}\,dy
\nonumber
\\
&&
-(u_{yy}\,q_{1,0}-(u_{xy}+\lambda)\,q_{0,1})\,dz,
\label{WEform}
\\
\omega_{i,j} &=& \widetilde{D}_x^i \widetilde{D}_y^j \,\omega_{0,0},
\qquad i,j \ge 0.
\nonumber
\end{eqnarray}
In this work we establish that the form $\omega_{0,0}$ can be derived from Maurer--Cartan forms of the contact symmetry pseudo-group of Eq. (\ref{PlebanskiEquation}).

\section{Symmetry pseudo-groups of differential equations}

Let $\pi :\mathbb{R}^n \times \mathbb{R} \rightarrow \mathbb{R}^n$ be a vector bundle with the local base coordinates $(x^1,...,x^n)$ and the local fibre coordinate $u$; then by $J^2(\pi)$ denote the bundle of the second-order jets of sections of $\pi$, with the local coordinates $(x^i,u,u_i,u_{ij})$, $i,j\in\{1,...,n\}$. For every local section $(x^i,f(x))$ of $\pi$, denote by $j_2(f)$ the corresponding 2-jet 
$(x^i,f(x),\partial f(x)/\partial x^i,\partial^2 f(x)/\partial x^i\partial x^j)$. A differential 1-form $\vartheta$ on $J^2(\pi)$ is called a {\it contact form} if it is annihilated by all 2-jets of local sections: $j_2(f)^{*}\vartheta = 0$. In the local coordinates every contact 1-form is a linear combination of the forms  $\vartheta_0 = du - u_{i}\,dx^i$,
$\vartheta_i = du_i - u_{ij}\,dx^j$, $i, j \in \{1,...,n\}$, $u_{ji} = u_{ij}$ (here and later we use the Einstein summation convention, so $u_i\,dx^i = \sum_{i=1}^{n}\,u_i\,dx^i$, etc.) A local diffeomorphism $\Delta : J^2(\pi) \rightarrow J^2(\pi)$, 
$\Delta : (x^i,u,u_i,u_{ij}) \mapsto (\bar{x}^i,\bar{u},\bar{u}_i,\bar{u}_{ij})$,
is called a {\it contact transformation} if for every contact 1-form $\bar{\vartheta}$ the form $\Delta^{*}\bar{\vartheta}$ is also contact. We denote by $\mathrm{Cont}(J^2(\pi))$  the pseudo-group  of contact transformations on $J^2(\pi)$.

Let $\EuScript{H} \subset \mathbb{R}^{(2 n+1)(n+3)(n+1)/3}$ be an open set with local coordinates
$a$, $b^i_k$, $c^i$, $f^{ik}$, $g_i$, $s_{ij}$, $w^k_{ij}$, $z_{ijk}$, $i,j,k \in \{1,...,n\}$,  such that $a\not =0$, $\det (b^i_k) \not = 0$, $f^{ik}=f^{ki}$, 
$z_{ijk}=z_{jik}= z_{ikj}$. Let $(B^i_k)$ be the inverse matrix for the matrix $(b^k_l)$, so
$B^i_k\,b^k_l = \delta^i_l$. We consider the {\it lifted coframe}
\[
\Theta_0 = a\, \vartheta_0,
\quad
\Theta_i = g_i\,\Theta_0 + a\,B_i^k\,\vartheta_k,
\quad
\Xi^i =c^i\,\Theta_0+f^{ik}\,\Theta_k+b_k^i\,dx^k,
\]
\begin{equation}
\Sigma_{ij} = s_{ij}\,\Theta_0+w_{ij}^{k}\,\Theta_k+z_{ijk}\,\Xi^k + a\,B^k_i\, B^l_j\,du_{kl},
\label{LCF}
\end{equation}
\noindent
$i \le j$, defined on $J^2(\pi)\times\EuScript{H}$. As it is shown in \cite{Morozov2006}, the forms (\ref{LCF}) are Maurer--Cartan forms for $\mathrm{Cont}(J^2(\pi))$, that is, a local diffeomorphism
$\widehat{\Delta} : J^2(\pi) \times \EuScript{H} \rightarrow J^2(\pi) \times \EuScript{H}$
satisfies the conditions
$\widehat{\Delta}^{*}\, \bar{\Theta}_0 = \Theta_0$,
$\widehat{\Delta}^{*}\, \bar{\Theta}_i = \Theta_i$,
$\widehat{\Delta}^{*}\, \bar{\Xi}^i = \Xi^i$,
and $\widehat{\Delta}^{*}\, \bar{\Sigma}_{ij} = \Sigma_{ij}$
whenever it is projectable on $J^2(\pi)$, and its projection
$\Delta : J^2(\pi) \rightarrow J^2(\pi)$ is a contact transformation.

The structure equations for $\mathrm{Cont}(J^2(\pi))$ read
\begin{eqnarray}
d \Theta_0 &=& \Phi^0_0 \wedge \Theta_0 + \Xi^i \wedge \Theta_i,
\nonumber
\\
d \Theta_i &=& \Phi^0_i \wedge \Theta_0 + \Phi^k_i \wedge \Theta_k
+ \Xi^k \wedge \Sigma_{ik},
\nonumber
\\
d \Xi^i &=& \Phi^0_0 \wedge \Xi^i -\Phi^i_k \wedge \Xi^k
+\Psi^{i0} \wedge \Theta_0
+\Psi^{ik} \wedge \Theta_k,
\nonumber
\\
d \Sigma_{ij} &=& \Phi^k_i \wedge \Sigma_{kj} - \Phi^0_0 \wedge \Sigma_{ij}
+ \Upsilon^0_{ij} \wedge \Theta_0
+ \Upsilon^k_{ij} \wedge \Theta_k + \Lambda_{ijk} \wedge \Xi^k,
\nonumber
\end{eqnarray}
\noindent
where the additional forms $\Phi^0_0$, $\Phi^0_i$, $\Phi^k_i$, $\Psi^{i0}$, $\Psi^{ij}$,
$\Upsilon^0_{ij}$, $\Upsilon^k_{ij}$, and $\Lambda_{ijk}$ depend on differentials of the coordinates of $\EuScript{H}$.

Suppose $\EuScript{E}$ is a second-order differential equation in one dependent and $n$ independent variables. We consider $\EuScript{E}$ as a sub-bundle in $J^2(\pi)$. Let $\mathrm{Cont}(\EuScript{E})$ be the group of contact symmetries for $\EuScript{E}$. It consists of all the contact transformations on $J^2(\pi)$ mapping $\EuScript{E}$ to itself. 
Let $\iota_0 : \EuScript{E} \rightarrow J^2(\pi)$ be an embedding and 
$\iota = \iota_0 \times \mathrm{id} : 
\EuScript{E}\times \EuScript{H} \rightarrow J^2(\pi)\times \EuScript{H}$. Maurer--Cartan forms of the pseudo-group $\mathrm{Cont}(\EuScript{E})$ can be obtained from the forms  $\theta_0 = \iota^{*} \Theta_0$, $\theta_i= \iota^{*}\Theta_i$, $\xi^i = \iota^{*}\Xi^i$ and $\sigma_{ij}=\iota^{*}\Sigma_{ij}$
by means of Cartan's method of equivalence, \cite{Cartan1}--\cite{Cartan4}, \cite{Gardner,Kamran,Olver95}, see details and examples in \cite{FelsOlver,Morozov2007,Morozov2008,Morozov2009}.

\section{Symmetry pseudo-group and the covering of Pleba\~nski's equation}

Following the method outlined in the previous section we find the Maurer--Cartan forms and their structure equations for the symmetry pseudo-group of Eq. (\ref{PlebanskiEquation}). The structure equations for the forms $\theta_0$, $\theta_i$,  $\xi^i$, $i \in \{1,2,3\}$, read  
\begin{eqnarray*}
%------------------------------------------------------------------------
\fl
d\theta_0 &=& \eta_5 \wedge \theta_0
+\xi^1 \wedge \theta_1
+\xi^2 \wedge \theta_2
+\xi^3 \wedge \theta_3
+\xi^4 \wedge \theta_4,
\\
%------------------------------------------------------------------------
\fl
d\theta_1&=&
(\eta_5-\eta_1) \wedge \theta_1
-\eta_3 \wedge \theta_4
-\eta_6 \wedge \theta_2
+\textfrac{1}{3}\,(\eta_5 
-4\,\eta_1 
+2\,\eta_4 
+3\,\sigma_{22}) \wedge \theta_3
\\
\fl
&&
+\xi^1 \wedge \sigma_{11}
+\xi^2 \wedge \sigma_{12}
+\xi^3 \wedge \sigma_{13}
+\xi^4 \wedge \sigma_{14},
\\
%------------------------------------------------------------------------
\fl
d\theta_2 &=& 
\textfrac{1}{3}\,(\eta_4 -2\,\eta_1 +2\,\eta_5) \wedge \theta_2
-\eta_3 \wedge \theta_3
+\xi^1 \wedge \sigma_{12}
+\xi^2 \wedge \sigma_{22}
+\xi^3 \wedge \sigma_{23}
\\
\fl
&&
+\xi^4 \wedge (\sigma_{13}+\sigma_{22}+\sigma_{33}),
\\
%------------------------------------------------------------------------
\fl
d\theta_3 &=& 
\textfrac{1}{3}\,(\eta_1 -2\,\eta_4 + 2\,\eta_5) \wedge \theta_3
-\eta_2 \wedge \theta_2 
+\xi^1 \wedge \sigma_{13}
+\xi^2 \wedge \sigma_{23}
+\xi^3 \wedge \sigma_{33}
+\xi^4 \wedge \sigma_{34},
\\
%------------------------------------------------------------------------
\fl
d\theta_4 &=& 
(\eta_5 -\eta_4) \wedge \theta_4 
-\eta_2 \wedge \theta_1
+\textfrac{1}{3}\,(\eta_5 -2\,\eta_1 +4\,\eta_4+ 3\,\sigma_{33}) \wedge \theta_2
+\xi^1 \wedge \sigma_{14}
\\
\fl
&&
-(2\,\eta_2+2\,\eta_3+\eta_6 - 2\,\sigma_{23}) \wedge \theta_3
+\xi^2 \wedge (\sigma_{13}+ \sigma_{22}+\sigma_{33})
+\xi^3 \wedge \sigma_{34}
+\xi^4 \wedge \sigma_{44},
\\
%------------------------------------------------------------------------
\fl
d\xi^1 &=& \eta_1 \wedge \xi^1+\eta_2 \wedge \xi^4,
\\
%------------------------------------------------------------------------
\fl
d\xi^2 &=& 
\eta_6 \wedge \xi^1
+\textfrac{1}{3}\,(\eta_5 + 2\,\eta_1 -\eta_4) \wedge \xi^2
+\eta_2 \wedge \xi^3
+\textfrac{1}{3}\,(\eta_5 -4\,\eta_4+2\,\eta_1+3\,\sigma_{33}) \wedge \xi^4,
\\
%------------------------------------------------------------------------
\fl
d\xi^3 &=&
\textfrac{1}{3}\,(4\,\eta_1-2\,\eta_4  -\eta_5-3\,\sigma_{22}) \wedge \xi^1
+\eta_3 \wedge \xi^2
+\textfrac{1}{3}\,(\eta_5  +2\,\eta_4-\eta_1) \wedge \xi^3
\\
\fl
&&
+(2\,\eta_2+2\,\eta_3+\eta_6 -2\,\sigma_{23}) \wedge \xi^4,
\\
%------------------------------------------------------------------------
\fl
d\xi^4 &=&\eta_3 \wedge \xi^1+\eta_4 \wedge \xi^4.
\end{eqnarray*}
The involutive system of structure equations for this pseudo-group is given in Appendix. In the next calculations  we use the following Maurer--Cartan forms only:
\begin{eqnarray*}
\xi^1 &=& b_{11}\,dt+b_{14}\,dz,
\\
\xi^2 &=& v^{-1}\,\left(
b_{11}\,dx+b_{14}\,dy
-(b_{11}\,(w-1)\,u_{xy}+b_{14}\,u_{xx}+b_{41}\,v)\,dt
\right.
\\
&&
\left.
-(b_{14}\,(w+1)\,u_{xy}-b_{11}\,u_{yy}+b_{44}\,v)\,dz
\right),
\\
\xi^3 &=& v^{-1}\,\left(
b_{41}\,dx+b_{44}\,dy
+(b_{11}\,v-b_{41}\,(w-1)\,u_{xy}-b_{44}\,u_{xx})\,dt
\right.
\\
&&
\left.
+(b_{14}\,v-b_{44}\,(w+1)\,u_{xy}+b_{41}\,u_{yy})\,dz
\right),
\\
\xi^4 &=& b_{41}\,dt+b_{44}\,dz,
\\
\eta_1 &=& 
(b_{44}\,db_{11}-b_{41}\,db_{14})\,(b_{11}b_{44}-b_{14}b_{41})^{-1}
+r_1\,\xi^1+r_2\,\xi^4,
\\
\eta_4 &=& (b_{11}\,db_{44}-b_{14}\,db_{41})\,(b_{11}b_{44}-b_{14}b_{41})^{-1}
-r_1\,\xi^1-r_2\,\xi^4,
\\
\eta_5 &=& -3\,v^{-1}\,dv+\eta_1+\eta_4,
\end{eqnarray*}
where $b_{11}$, $b_{14}$, $b_{41}$, $b_{44}$, $v$, $w$, $r_1$, $r_2$ are arbitrary parameters such that  $b_{11}b_{44}-b_{14}b_{41}\not =0$ and  $v\not =0$.

Direct computations prove the following

\vskip 10 pt
\noindent
{\sc theorem.}
{\it
Either substituting for  $v=q_{0,0}$, $b_{11}=q_{1,0}$, $b_{14}=q_{0,1}$, $w=\lambda\,u_{xy}^{-1}$
into the linear combination
\[
\textfrac{1}{3}\,(\eta_1+\eta_4 -\eta_5)-\xi^2-\xi^4,
\]
or substituting for $v=q_{0,0}$, $b_{41}=q_{1,0}$, $b_{44}=q_{0,1}$, $w=\lambda\,u_{xy}^{-1}$
into the linear combination
\[
\textfrac{1}{3}\,(\eta_1+\eta_4 -\eta_5)+\xi^1-\xi^3
\]
yields the form $q_{0,0}^{-1}\,\omega_{0,0}$ proportional to the form (\ref{WEform}), which is the generating form of the ideal of Wahquist--Estabrook forms of the covering (\ref{covering}) of Eq. (\ref{PlebanskiEquation}).
}

\vskip 10 pt

Another approach to computing Wahquist--Estabrook forms of coverings of {\sc pde}s from Maurer--Cartan forms of their symmetry pseudo-groups was proposed in \cite{Morozov2008b}. We hope to apply this to Eq. (\ref{PlebanskiEquation}) elsewhere.

\vskip 20 pt
%\newpage %%%%%%%%%%%%%%%%%%%%%%%%%%%%%%%%%%%%%%%%%%%%%%%%%%%%%%%%%%%%%%%%%%%%%%%%%%%%%%%%%%%%%
\section*{Appendix}

The involutive system of structure equations for the symmetry pseudo-group of Eq. (\ref{PlebanskiEquation}):
\begin{eqnarray}
%------------------------------------------------------------------------
\fl
d\theta_0 &=& \eta_5 \wedge \theta_0
+\xi^1 \wedge \theta_1
+\xi^2 \wedge \theta_2
+\xi^3 \wedge \theta_3
+\xi^4 \wedge \theta_4,
\nonumber \\
%------------------------------------------------------------------------
\fl
d\theta_1&=&
(\eta_5-\eta_1) \wedge \theta_1
-\eta_3 \wedge \theta_4
-\eta_6 \wedge \theta_2
+\textfrac{1}{3}\,(\eta_5 
-4\,\eta_1 
+2\,\eta_4 
+3\,\sigma_{22}) \wedge \theta_3
\nonumber \\
\fl
&&
+\xi^1 \wedge \sigma_{11}
+\xi^2 \wedge \sigma_{12}
+\xi^3 \wedge \sigma_{13}
+\xi^4 \wedge \sigma_{14},
\nonumber \\
%------------------------------------------------------------------------
\fl
d\theta_2 &=& 
\textfrac{1}{3}\,(\eta_4 -2\,\eta_1 +2\,\eta_5) \wedge \theta_2
-\eta_3 \wedge \theta_3
+\xi^1 \wedge \sigma_{12}
+\xi^2 \wedge \sigma_{22}
+\xi^3 \wedge \sigma_{23}
\nonumber \\
\fl
&&
+\xi^4 \wedge (\sigma_{13}+\sigma_{22}+\sigma_{33}),
\nonumber \\
%------------------------------------------------------------------------
\fl
d\theta_3 &=& 
\textfrac{1}{3}\,(\eta_1 -2\,\eta_4 + 2\,\eta_5) \wedge \theta_3
-\eta_2 \wedge \theta_2 
+\xi^1 \wedge \sigma_{13}
+\xi^2 \wedge \sigma_{23}
+\xi^3 \wedge \sigma_{33}
+\xi^4 \wedge \sigma_{34},
\nonumber \\
%------------------------------------------------------------------------
\fl
d\theta_4 &=& 
(\eta_5 -\eta_4) \wedge \theta_4 
-\eta_2 \wedge \theta_1
+\textfrac{1}{3}\,(\eta_5 -2\,\eta_1 +4\,\eta_4+ 3\,\sigma_{33}) \wedge \theta_2
+\xi^1 \wedge \sigma_{14}
\nonumber \\
\fl
&&
-(2\,\eta_2+2\,\eta_3+\eta_6 - 2\,\sigma_{23}) \wedge \theta_3
+\xi^2 \wedge (\sigma_{13}+ \sigma_{22}+\sigma_{33})
+\xi^3 \wedge \sigma_{34}
+\xi^4 \wedge \sigma_{44},
\nonumber \\
%------------------------------------------------------------------------
\fl
d\xi^1 &=& \eta_1 \wedge \xi^1+\eta_2 \wedge \xi^4,
\nonumber \\
%------------------------------------------------------------------------
\fl
d\xi^2 &=& 
\eta_6 \wedge \xi^1
+\textfrac{1}{3}\,(\eta_5 + 2\,\eta_1 -\eta_4) \wedge \xi^2
+\eta_2 \wedge \xi^3
+\textfrac{1}{3}\,(\eta_5 -4\,\eta_4+2\,\eta_1+3\,\sigma_{33}) \wedge \xi^4,
\nonumber \\
%------------------------------------------------------------------------
\fl
d\xi^3 &=&
\textfrac{1}{3}\,(4\,\eta_1-2\,\eta_4  -\eta_5-3\,\sigma_{22}) \wedge \xi^1
+\eta_3 \wedge \xi^2
+\textfrac{1}{3}\,(\eta_5  +2\,\eta_4-\eta_1) \wedge \xi^3
\nonumber \\
\fl
&&
+(2\,\eta_2+2\,\eta_3+\eta_6 -2\,\sigma_{23}) \wedge \xi^4,
\nonumber \\
%------------------------------------------------------------------------
\fl
d\xi^4 &=&\eta_3 \wedge \xi^1+\eta_4 \wedge \xi^4,
\nonumber \\
%------------------------------------------------------------------------
\fl
d\sigma_{11}&=&
\eta_7 \wedge \theta_1 
+\eta_8 \wedge \theta_2
+\eta_9 \wedge \theta_3
+\eta_{10} \wedge \theta_4
+\eta_{11} \wedge \xi^1
+\eta_{12} \wedge \xi^2
+\eta_{13} \wedge \xi^3
+\eta_{14} \wedge \xi^4 
\nonumber \\
\fl
&&
+(\eta_5 -2\,\eta_1) \wedge \sigma_{11}
-2\,\eta_6 \wedge \sigma_{12}
+\textfrac{2}{3}\,(\eta_5 -4\,\eta_1+2\,\eta_4 +3\,\sigma_{22}) \wedge \sigma_{13}
-2\,\eta_3 \wedge \sigma_{14},
\nonumber \\
%------------------------------------------------------------------------
\fl
d\sigma_{12}&=&
\eta_7 \wedge \theta_2
+\eta_{10} \wedge \theta_3
+\eta_{12} \wedge \xi^1
+(2\,\eta_7 -\eta_9) \wedge \xi^2
+\eta_{15} \wedge \xi^3
+\eta_{16} \wedge \xi^4
-(\eta_3 +\eta_6) \wedge \sigma_{22}
\nonumber \\
\fl
&&
+\textfrac{1}{3}\,(2\,\eta_5+\eta_4-5\,\eta_1) \wedge \sigma_{12}
-2\,\eta_3 \wedge \sigma_{13}
+\textfrac{1}{3}\,(\eta_5 +2\,\eta_4-4\,\eta_1+3\,\sigma_{2\,2}) \wedge \sigma_{23}
-\eta_3 \wedge \sigma_{33},
\nonumber \\
%------------------------------------------------------------------------
\fl
d\sigma_{13} &=&
(\eta_{15} -\eta_8 -\eta_{10}) \wedge \theta_2
-\eta_7 \wedge \theta_3 
+\eta_{13} \wedge \xi^1
+\eta_{15} \wedge \xi^2
+(\eta_{16} -2\,\eta_7+\eta_9-\eta_{13}) \wedge \xi^3
\nonumber \\
\fl
&&
+\eta_{17} \wedge \xi^4
-\eta_2 \wedge \sigma_{12}
+\textfrac{2}{3}\,(\eta_5-\eta_1-\eta_4) \wedge \sigma_{13}
-\eta_6 \wedge \sigma_{23} 
-\eta_3 \wedge \sigma_{34}
\nonumber \\
\fl
&&
+\textfrac{1}{3}\,(\eta_5 -4\,\eta_1+2\,\eta_4+3\,\sigma_{22}) \wedge \sigma_{33},
\nonumber \\
%------------------------------------------------------------------------
\fl
d\sigma_{14} &=&
(\eta_{15}-\eta_8-\eta_{10}) \wedge \theta_1
+(\eta_{16}+\eta_9-\eta_{13}) \wedge \theta_2
-\eta_8 \wedge \theta_3
-\eta_7 \wedge \theta_4
+\eta_{14} \wedge \xi^1
+\eta_{16} \wedge \xi^2
\nonumber \\
\fl
&&
+\eta_{17} \wedge \xi^3
+\eta_{18} \wedge \xi^4
-\eta_2 \wedge \sigma_{11}
-\textfrac{1}{3}\,(\eta_5 -2\,\eta_1+4\,\eta_4 -3\,\sigma_{33}) \wedge \sigma_{12}
\nonumber \\
\fl
&&
-2\,(\eta_2+\eta_3+\eta_6-\sigma_{23}) \wedge \sigma_{13}
+(\eta_5-\eta_1-\eta_4) \wedge \sigma_{14}
-\eta_6 \wedge (\sigma_{22}+\sigma_{33})
\nonumber \\
\fl
&&
+\textfrac{1}{3}\,(\eta_5 -4\,\eta_1 +2\,\eta_4+3\,\sigma_{22}) \wedge \sigma_{34}
-\eta_3 \wedge \sigma_{44},
\nonumber \\
%------------------------------------------------------------------------
\fl
d\sigma_{22} &=& 
(2\,\eta_7-\eta_9) \wedge \xi^1
-\eta_{10} \wedge \xi^2
+\eta_7 \wedge \xi^3
+(2\,\eta_{15}-2\,\eta_{10}-\eta_8) \wedge \xi^4
-2\,\eta_3 \wedge \sigma_{23}
\nonumber \\
\fl
&&
+\textfrac{1}{3}\,(\eta_5 -4\,\eta_1+2\,\eta_4) \wedge \sigma_{22},
\nonumber \\
%------------------------------------------------------------------------
\fl
d\sigma_{23} &=& 
\eta_{15} \wedge \xi^1
+\eta_7 \wedge \xi^2 
+(\eta_{15}-\eta_8-\eta_{10}) \wedge \xi^3
+\eta_{19} \wedge \xi^4
-\eta_2 \wedge \sigma_{22}
-\eta_3 \wedge \sigma_{33}
\nonumber \\
\fl
&&
+\textfrac{1}{3}\,(\eta_5-\eta_1 -\eta_4) \wedge \sigma_{23},
\nonumber \\
%------------------------------------------------------------------------
\fl
d\sigma_{33} &=& 
(\eta_{16}-2\,\eta_7+\eta_9-\eta_{13}) \wedge \xi^1
+(\eta_{15}-\eta_8-\eta_{10}) \wedge \xi^2
+(\eta_{19}+\eta_7-\eta_9+\eta_{13}-\eta_{16}) \wedge \xi^3
\nonumber \\
\fl
&&
+\eta_{20} \wedge \xi^4
-2\,\eta_2 \wedge \sigma_{23}
+\textfrac{1}{3}\,(\eta_5+2\,\eta_1-4\,\eta_4) \wedge \sigma_{33},
\nonumber \\
%------------------------------------------------------------------------
\fl
d\sigma_{34} &=& 
(\eta_7-\eta_9+\eta_{13}-\eta_{16}+\eta_{19}) \wedge \theta_2 
+(\eta_8 +\eta_{10} -\eta_{15}) \wedge \theta_3
+\eta_{17} \wedge \xi^1
+\eta_{19} \wedge \xi^2
+\eta_{20} \wedge \xi^3
\nonumber \\
\fl
&&
+\eta_{21} \wedge \xi^4
-2\,\eta_2 \wedge \sigma_{13}
-\eta_2 \wedge \sigma_{22}
-\textfrac{1}{3}\,(\eta_5 -2\,\eta_1+4\,\eta_4) \wedge \sigma_{23}
\nonumber \\
\fl
&&
-(3\,\eta_2 +2\,\eta_3+\eta_6-3\,\sigma_{23}) \wedge \sigma_{33}
+\textfrac{1}{3}\,(2\,\eta_5+\eta_1-5\,\eta_4) \wedge \sigma_{34},
\nonumber %
\\
%\end{eqnarray}
%%%%%%%%%%%%%%%%%%%%%%%%%%%%%%%%%%%%%%%%%%%%%%%%%%%%%%%%%%%%%%%%%%%%%%%%%%%%%%%%%%%%
%\begin{eqnarray}
%------------------------------------------------------------------------
\fl
d\sigma_{44} &=& 
(\eta_7 -\eta_9 +\eta_{13}-\eta_{16}+\eta_{19}) \wedge \theta_1
+(\eta_{20}-2\,(\eta_8+\eta_{10}-\eta_{15})) \wedge \theta_2
+\eta_{18} \wedge \xi^1 
\nonumber \\
\fl
&&
-(\eta_9 -\eta_{13}+\eta_{16}) \wedge \theta_3
+(\eta_8 +\eta_{10}-\eta_{15}) \wedge \theta_4
+(\eta_{20}-\eta_8-2\,\eta_{10}+2\,\eta_{15}+\eta_{17}) \wedge \xi^2
\nonumber \\
\fl
&&
+\eta_{21} \wedge \xi^3
+\eta_{22} \wedge \xi^4
-\textfrac{2}{3}\,(\eta_5 +2\,\eta_1-4\,\eta_4+6\,\sigma_{33}) \wedge (\sigma_{13}+\sigma_{33})
-2\,\eta_2 \wedge \sigma_{14}
\nonumber \\
\fl
&&
-\textfrac{2}{3}\,(\eta_5 +2\,\eta_1-4\,\eta_4+6\,\sigma_{33}) \wedge \sigma_{22}
-2\,(\eta_6+2\,(\eta_2+\eta_3-\sigma_{23})) \wedge \sigma_{34}
\nonumber \\
\fl
&&
+(\eta_5-2\,\eta_4) \wedge \sigma_{44},
\nonumber \\
%------------------------------------------------------------------------
\fl
d\eta_1 &=&  
\eta_7 \wedge \xi^1
+(\eta_{15}-\eta_8-\eta_{10}) \wedge \xi^4
+\eta_2 \wedge \eta_3,
\nonumber \\
%------------------------------------------------------------------------
\fl
d\eta_2 &=& 
(\eta_{15}-\eta_8-\eta_{10}) \wedge \xi^1
+(\eta_{19}+\eta_7-\eta_9+\eta_{13}-\eta_{16}) \wedge \xi^4
+(\eta_1 -\eta_4) \wedge \eta_2,
\nonumber \\
%------------------------------------------------------------------------
\fl
d\eta_3 &=& 
\eta_{10} \wedge \xi^1
-\eta_7 \wedge \xi^4
+(\eta_4-\eta_1) \wedge \eta_3,
\nonumber \\
%------------------------------------------------------------------------
\fl
d\eta_4 &=& 
(\eta_8 +\eta_{10}-\eta_{15}) \wedge \xi^4
-\eta_7 \wedge \xi^1
-\eta_2 \wedge \eta_3,
\nonumber \\
%------------------------------------------------------------------------
\fl
d\eta_5 &=& 0,
\nonumber \\
%------------------------------------------------------------------------
\fl
d\eta_6 &=&
\eta_8 \wedge \xi^1
+(\eta_7-\eta_8-\eta_{10}+\eta_{15}) \wedge \xi^3
+(\eta_9 -\eta_{13}+\eta_{16}) \wedge \xi^4
+\textfrac{1}{3}\,(\eta_6 - 4\,\eta_2 -2\,\eta_3) \wedge \eta_1
\nonumber \\
\fl
&&
+\textfrac{1}{3}\,(2\,\eta_4+\eta_5+3\,\sigma_{22}) \wedge \eta_2 
+\textfrac{1}{3}\,(\eta_5 -4\,\eta_4+3\,\sigma_{33}) \wedge \eta_3
-\textfrac{1}{3}(\eta_4-\eta_5) \wedge \eta_6,
\nonumber \\
%------------------------------------------------------------------------
\fl
d\eta_7 &=& 
\eta_{23} \wedge \xi^1
+\eta_{24} \wedge \xi^4
+\eta_7 \wedge \eta_1
-\eta_{10} \wedge \eta_2
+2\,(\eta_{15}-\eta_8-\eta_{10}) \wedge \eta_3,
\nonumber \\
%------------------------------------------------------------------------
\fl
d\eta_8 &=&  
\eta_{25} \wedge \xi^1 
+\eta_{23} \wedge \xi^2
+\eta_{24} \wedge \xi^3
+\eta_{26} \wedge \xi^4
-\textfrac{2}{3}\,(2\,\eta_8 + 5\,\eta_{10}-4\,\eta_{15}) \wedge \eta_1
+\eta_2 \wedge \eta_9
\nonumber \\
\fl
&&
+2\,(\eta_9-\eta_{13}+\eta_{16}) \wedge \eta_3
+\textfrac{1}{3}\,(5\,\eta_8+8\,\eta_{10}-4\,\eta_{15}) \wedge \eta_4
+\textfrac{1}{3}\,(\eta_8 +\eta_{10}-2\,\eta_{15}) \wedge \eta_5
\nonumber \\
\fl
&&
-\eta_6 \wedge \eta_7
+ 2\,(\eta_8+\eta_{10}-\eta_{15}) \wedge \sigma_{22}
-\eta_{10} \wedge \sigma_{33},
\nonumber \\
%------------------------------------------------------------------------
\fl
d\eta_9 &=& 
\eta_{27} \wedge \xi^1
+\eta_{28} \wedge \xi^2
-\eta_{23} \wedge \xi^3
-\eta_{25} \wedge \xi^4
-\textfrac{1}{3}\,(12\,\eta_7 +7\,\eta_9) \wedge \eta_1
-(3\,\eta_8+2\,\eta_{10}) \wedge \eta_3
\nonumber \\
\fl
&&
+\eta_{10} \wedge (\eta_6+2\,\sigma_{23}-2\,\eta_2)
\nonumber \\
\fl
&&
+\textfrac{2}{3}\,(3\,\eta_7-\eta_9) \wedge \eta_4
+\textfrac{1}{3}\,(3\,\eta_7 -\eta_9) \wedge \eta_5
+3\,\eta_7 \wedge \sigma_{22},
\nonumber \\
%------------------------------------------------------------------------
\fl
d\eta_{10} &=& 
\eta_{28} \wedge \xi^1
-\eta_{23} \wedge \xi^4
+\eta_{10} \wedge (2\,\eta_1-\eta_4)
+3\,\eta_3 \wedge \eta_7,
\nonumber \\
%------------------------------------------------------------------------
\fl
d\eta_{11} &=&
\eta_{23} \wedge \theta_1
+\eta_{25} \wedge \theta_2
+\eta_{27} \wedge \theta_3
+\eta_{28} \wedge \theta_4
+\eta_{29} \wedge \xi^1
+\eta_{30} \wedge \xi^2
+\eta_{31} \wedge \xi^3
+\eta_{32} \wedge \xi^4
\nonumber \\
\fl
&&
+(3\,\eta_{11}+4\,\eta_{13}) \wedge \eta_1
+3\,\eta_{14} \wedge \eta_3
-2\,\eta_{13} \wedge \eta_4
-(\eta_{11}+\eta_{13}) \wedge \eta_5
\nonumber \\
\fl
&&
+3\,(\eta_{12} \wedge (\eta_6+\eta_8)
-\eta_7 \wedge \sigma_{11}
-\eta_9 \wedge \sigma_{13}
-\eta_{10} \wedge \sigma_{14} 
-\eta_{13} \wedge \sigma_{22}),
\nonumber \\
%------------------------------------------------------------------------
\fl
d\eta_{12} &=& 
\eta_{23} \wedge \theta_2
+\eta_{28} \wedge \theta_3
+\eta_{30} \wedge \xi^1 
+(2\,\eta_{23}-\eta_{27}) \wedge \xi^2
+(\eta_{24}+\eta_{25}+\eta_{28}) \wedge \xi^3
\nonumber \\
\fl
&&
+(\eta_{26}-\eta_{27}+\eta_{31}) \wedge \xi^4
+\textfrac{8}{3}\,(\eta_{12}+\eta_{15}) \wedge \eta_1
+(\eta_{13}+2\,\eta_{16}) \wedge \eta_3
-\textfrac{1}{3}\,(\eta_{12} +4\,\eta_{15}) \wedge \eta_4
\nonumber \\
\fl
&&
-\textfrac{2}{3}\,(\eta_{12}+\eta_{15}) \wedge \eta_5
+2\,(2\,\eta_7 -\eta_9) \wedge \eta_6
-3\,\eta_7 \wedge \sigma_{12}
-(\eta_8+\eta_{10}+2\,\eta_{15}) \wedge \sigma_{22}
\nonumber \\
\fl
&&
-\eta_9 \wedge \sigma_{23}
-\eta_{10} \wedge (3\,\sigma_{13}+\sigma_{33}),
\nonumber \\
%------------------------------------------------------------------------
\fl
d\eta_{13} &=& 
\eta_{24} \wedge \theta_2
-\eta_{23} \wedge \theta_3
+\eta_{31} \wedge \xi^1
+(\eta_{24}+\eta_{25}+\eta_{28}) \wedge \xi^2
-(2\,\eta_{23}-\eta_{26}) \wedge \xi^3
+\eta_{33} \wedge \xi^4
\nonumber \\
\fl
&&
-\textfrac{1}{3}\,(16\,\eta_7+8\,\eta_9-3\,\eta_{13}+8\,\eta_{16}) \wedge \eta_1
+\eta_{12} \wedge \eta_2
+2\,\eta_{17} \wedge \eta_3
\nonumber \\
\fl
&&
+\textfrac{2}{3}\,(4\,\eta_7 -2\,\eta_9 +3\,\eta_{13}-2\,\eta_{16}) \wedge \eta_4
+\textfrac{2}{3}\,(2\,\eta_7-\eta_9-\eta_{16}) \wedge \eta_5
+2\,\eta_{15} \wedge(\eta_6-\sigma_{12})
\nonumber \\
\fl
&&
+\eta_7 \wedge (\sigma_{13}+4\,\sigma_{22})
+\eta_8 \wedge (2\,\sigma_{12}-\sigma_{23})
-\eta_9 \wedge (2\, \sigma_{22}+\sigma_{33})
+\eta_{10} \wedge (2\,\sigma_{12}-\sigma_{34})
\nonumber \\
\fl
&&
+2\,(\eta_{13} -\eta_{16}) \wedge \sigma_{22},
\nonumber \\
%------------------------------------------------------------------------
\fl
d\eta_{14} &=&
\eta_{24} \wedge \theta_1
+\eta_{26} \wedge \theta_2
-\eta_{25} \wedge \theta_3
-\eta_{23} \wedge \theta_4
+\eta_{32} \wedge \xi^1
+(\eta_{26}-\eta_{27}+\eta_{31}) \wedge \xi^2
+\eta_{33} \wedge \xi^3
\nonumber \\
\fl
&&
+\eta_{34} \wedge \xi^4
+\textfrac{2}{3}\,(\eta_{12}+3\,\eta_{14}+4\,\eta_{17}) \wedge \eta_1
+(\eta_{11}+2\,\eta_{13}) \wedge \eta_2
+2\,(\eta_{13}+\eta_{18}) \wedge \eta_3
\nonumber \\
\fl
&&
-\textfrac{1}{3}\,(4\,\eta_{12}-3\,\eta_{14}+4\,\eta_{17}) \wedge \eta_4
+\textfrac{1}{3}\,(\eta_{12}-3\,\eta_{14}-2\,\eta_{17}) \wedge \eta_5
+(\eta_{13}+2\,\eta_{16}) \wedge \eta_6
\nonumber \\
\fl
&&
+\eta_7 \wedge \sigma_{14}
+\eta_8 \wedge (2\,\sigma_{11}+\sigma_{13}-\sigma_{22}-\sigma_{33})
-\eta_9 \wedge (2\,\sigma_{12}+\sigma_{34})
+\eta_{10} \wedge (2\,\sigma_{11}-\sigma_{44})
\nonumber \\
\fl
&&
+\eta_{12} \wedge \sigma_{33}
+2\,((\eta_{13}+\eta_{16}) \wedge \sigma_{12}
-\eta_{13} \wedge \sigma_{23}
-\eta_{15} \wedge \sigma_{11}
-\eta_{17} \wedge \sigma_{22}),
\nonumber %
\\
%\end{eqnarray}
%%%%%%%%%%%%%%%%%%%%%%%%%%%%%%%%%%%%%%%%%%%%%%%%%%%%%%%%%%%%%%%%%%%%%%%%%%%%%%%%%%%%
%\begin{eqnarray}
%------------------------------------------------------------------------
\fl
d\eta_{15} &=& 
(\eta_{24}+\eta_{25}+\eta_{28}) \wedge \xi^1
+\eta_{23} \wedge \xi^2
+\eta_{24} \wedge \xi^3
+\eta_{35} \wedge \xi^4
-\textfrac{4}{3}\,(\eta_8 -\eta_{10}+2\,\eta_{15}) \wedge \eta_1
\nonumber \\
\fl
&&
+(2\,\eta_7-\eta_9) \wedge \eta_2
-(2\,\eta_7-\eta_9+\eta_{13}-\eta_{16}-\eta_{19}) \wedge \eta_3
+\textfrac{1}{3}\,(2\,(\eta_8+\eta_{10})-\eta_{15}) \wedge \eta_4
\nonumber \\
\fl
&&
+\textfrac{1}{3}\,(\eta_8 +\eta_{10}-2\,\eta_{15}) \wedge \eta_5
-\eta_6 \wedge \eta_7
+2\,(\eta_8 +\eta_{10}-\eta_{15}) \wedge \sigma_{22}
-\eta_{10} \wedge \sigma_{33},
\nonumber \\
%------------------------------------------------------------------------
\fl
d\eta_{16} &=&
\eta_{24} \wedge \theta_2
-\eta_{23} \wedge \theta_3
+(\eta_{26}-\eta_{27}+\eta_{31}) \wedge \xi^1
+(2\,\eta_{24}+\eta_{25}) \wedge \xi^2
+\eta_{35} \wedge \xi^3
+\eta_{36} \wedge \xi^4
\nonumber \\
\fl
&&
+\textfrac{1}{3}\,(4\,\eta_7-2\,\eta_9+5\,\eta_{16}+4\,\eta_{19}) \wedge \eta_1
+(\eta_{12}+2\,\eta_{15}) \wedge \eta_2
-(\sigma_{13}+2\,\sigma_{33}) \wedge \eta_7 
\nonumber \\
\fl
&&
-(\eta_8 +2\,\eta_{10} -4\,\eta_{15}-2\,\eta_{17}-\eta_{20}) \wedge \eta_3
-\textfrac{2}{3}\,(4\,\eta_7 -2\,\eta_9-\eta_{16}+\eta_{19}) \wedge \eta_4
\nonumber \\
\fl
&&
+\textfrac{1}{3}\,(2\,\eta_7-\eta_9-2\,\eta_{16}-\eta_{19}) \wedge \eta_5
-(\eta_8 +2\,\eta_{10}-3\,\eta_{15}) \wedge \eta_6
-(2\,\sigma_{12}+\sigma_{23}) \wedge \eta_8
\nonumber \\
\fl
&&
+(\sigma_{22}+\sigma_{33}) \wedge \eta_9
-(2\,\sigma_{12}+\sigma_{34}) \wedge \eta_{10}
+(\eta_{13}-\eta_{16}-\eta_{19}) \wedge \sigma_{22}
-2\,\eta_{15} (\wedge \sigma_{12}+\sigma_{23}),
\nonumber \\
%------------------------------------------------------------------------
\fl
d\eta_{17} &=& 
(\eta_{23}-\eta_{26}+\eta_{35}) \wedge \theta_2
-\eta_{24} \wedge \theta_3
+\eta_{33} \wedge \xi^1
+\eta_{35} \wedge \xi^2
-(2\,\eta_{24}+\eta_{25}+\eta_{33}-\eta_{36}) \wedge \xi^3
\nonumber \\
\fl
&&
+\eta_{37} \wedge \xi^4
+\textfrac{2}{3}\,(\eta_{15}+\eta_{17}+2\,\eta_{20}) \wedge \eta_1
-(4\,\eta_7-2\,\eta_9+\eta_{13}-3\,\eta_{16}) \wedge \eta_2
\nonumber \\
\fl
&&
-(4\,\eta_7-2\,\eta_9+2\,\eta_{13}-2\,\eta_{16}-\eta_{21}) \wedge \eta_3
-\textfrac{1}{3}\,(4\,\eta_{15}-5\,\eta_{17}+2\,\eta_{20}) \wedge \eta_4
\nonumber \\
\fl
&&
+\textfrac{1}{3}\,(\eta_{15}-2\,\eta_{17}-\eta_{20}) \wedge \eta_5
-(2\,\eta_7-\eta_9+\eta_{13}-\eta_{16}-\eta_{19}) \wedge \eta_6
-(\sigma_{12}+3\,\sigma_{23}) \wedge \eta_9
\nonumber \\
\fl
&&
+(\sigma_{12}-4\,\sigma_{23}-2\,\sigma_{34}) \wedge \eta_7
-(\sigma_{13}-\sigma_{22}-2\,\sigma_{33}) \wedge \eta_8
-(\sigma_{13}+\sigma_{22}+\sigma_{33}) \wedge \eta_{10}
\nonumber \\
\fl
&&
+(\sigma_{12}-3\,\sigma_{23}) \wedge \eta_{13} 
+(\sigma_{13}+\sigma_{22}) \wedge \eta_{15}
-(\sigma_{12}+3\,\sigma_{23}) \wedge \eta_{16}
-\eta_{19} \wedge \sigma_{12}
-\eta_{20} \wedge \sigma_{22},
\nonumber \\
%------------------------------------------------------------------------
\fl
d\eta_{18}  &=&  
(\eta_{23}-\eta_{26}+\eta_{35}) \wedge \theta_1
-(\eta_{25}+\eta_{33}-\eta_{36}) \wedge \theta_2
-\eta_{26} \wedge \theta_3
-\eta_{24} \wedge \theta_4
+\eta_{34} \wedge \xi^1
\nonumber \\
\fl
&&
+\eta_{36} \wedge \xi^2
+\eta_{37} \wedge \xi^3
+\eta_{38} \wedge \xi^4
+\textfrac{1}{3}\,(4\,\eta_{16}+3\,\eta_{18}+4\,\eta_{21}) \wedge \eta_1
+2\,(\eta_{14}+2\,\eta_{17}) \wedge \eta_2
\nonumber \\
\fl
&&
+(4\,\eta_{17}+\eta_{22}) \wedge \eta_3
-\textfrac{2}{3}\,(4\,\eta_{16}+3\,\eta_{18}-\eta_{21}) \wedge \eta_4
+\textfrac{1}{3}\,(2\,\eta_{16}-3\,\eta_{18}-\eta_{21}) \wedge \eta_5
\nonumber \\
\fl
&&
-(\eta_8 +2\,\eta_{10}-2\,\eta_{15}-3\,\eta_{17}-\eta_{20}) \wedge \eta_6
-\eta_7 \wedge (\sigma_{11}-2\,\sigma_{44})
+\eta_8 \wedge (2\,\sigma_{12}+\sigma_{14}+2\,\sigma_{34})
\nonumber \\
\fl
&&
+\eta_9 \wedge (\sigma_{11}-\sigma_{13}-2\,\sigma_{22}-2\,\sigma_{33})
+\eta_{10} \wedge (2\,\sigma_{12}+ \sigma_{14})
-\eta_{15} \wedge (2\,\sigma_{12}+\sigma_{14})
\nonumber \\
\fl
&&
-\eta_{13} \wedge (\sigma_{11}-\sigma_{13}-2\,\sigma_{22}-2\,\sigma_{33})
+\eta_{16} \wedge (\sigma_{11}-\sigma_{13}-2\,\sigma_{22})
-4\,\eta_{17} \wedge \sigma_{23}
-\eta_{19} \wedge \sigma_{11}
\nonumber \\
\fl
&&
-\eta_{20} \wedge \sigma_{12}
-\eta_{21} \wedge \sigma_{22},
\nonumber \\
%------------------------------------------------------------------------
\fl
d\eta_{19} &=& 
\eta_{35} \wedge \xi^1
+\eta_{24} \wedge \xi^2
+(\eta_{23}-\eta_{26}+\eta_{35}) \wedge \xi^3
+\eta_{39} \wedge \xi^4
+\textfrac{1}{3}\,(2\,\eta_7+\eta_{19}) \wedge \eta_1
\nonumber \\
\fl
&&
-(3\,\eta_8+4\,\eta_{10}-5\,\eta_{15}) \wedge \eta_2
-(2\,\eta_8+2\,\eta_{10}-2\,\eta_{15}-\eta_{20}) \wedge \eta_3
-\textfrac{4}{3}\,(\eta_7-\eta_{19}) \wedge \eta_4
\nonumber \\
\fl
&&
+\textfrac{1}{3}\,(\eta_7 -\eta_{19}) \wedge \eta_5
-(\eta_8+\eta_{10}-\eta_{15}) \wedge \eta_6
-\eta_7 \wedge (\sigma_{22}-2\,\sigma_{33})
\nonumber \\
\fl
&&
+2\,(\eta_8 +\eta_{10}-\eta_{15}) \wedge \sigma_{23}
+(\eta_9-\eta_{13}+\eta_{16}-\eta_{19}) \wedge \sigma_{22},
\nonumber \\
%------------------------------------------------------------------------
\fl
d\eta_{20} &=&   
(\eta_{36}-2\,\eta_{24}-\eta_{25}-\eta_{33}) \wedge \xi^1
+(\eta_{23}-\eta_{26}+\eta_{35}) \wedge \xi^2
-(\eta_9-\eta_{13}+\eta_{16}-4\,\eta_{19}) \wedge \eta_2
\nonumber \\
\fl
&&
+(\eta_{24}+\eta_{25}+\eta_{33}-\eta_{36}+\eta_{39}) \wedge \xi^3
+\eta_{40} \wedge \xi^4
-\textfrac{2}{3}\,(\eta_8+\eta_{10}-\eta_{15}+\eta_{20}) \wedge \eta_1
\nonumber \\
\fl
&&
+2\,(\eta_7-\eta_9+\eta_{13}-\eta_{16}+\eta_{19}) \wedge \eta_3
+\textfrac{1}{3}\,(4\,(\eta_8+\eta_{10}-\eta_{15})+7\,\eta_{20}) \wedge \eta_4
\nonumber \\
\fl
&&
-\textfrac{1}{3}\,(\eta_8+\eta_{10}-\eta_{15}+\eta_{20}) \wedge \eta_5
+(\eta_7-\eta_9+\eta_{13}-\eta_{16}+\eta_{19}) \wedge \eta_6
\nonumber \\
\fl
&&
-4\,(\eta_7 -\eta_9+\eta_{13}-\eta_{16}+\eta_{19}) \wedge \sigma_{23}
-3\,(\eta_8+\eta_{10}-\eta_{15}) \wedge \sigma_{33},
\nonumber \\
%------------------------------------------------------------------------
\fl
d\eta_{21} &=&
(\eta_{24}+\eta_{25}+\eta_{33}-\eta_{36}+\eta_{39}) \wedge \theta_2
-(\eta_{23}-\eta_{26}+\eta_{35}) \wedge \theta_3
+\eta_{37} \wedge \xi^1
+\eta_{39} \wedge \xi^2
\nonumber \\
\fl
&&
+\eta_{40} \wedge \xi^3
+\eta_{41} \wedge \xi^4
+\textfrac{1}{3}\,(4\,\eta_{19}-\eta_{21}) \wedge \eta_1
-(\eta_8+2\,\eta_{10}-2\,\eta_{15}-3\,\eta_{17}-5\,\eta_{20}) \wedge \eta_2
\nonumber \\
\fl
&&
-\textfrac{2}{3}\,(\eta_{19}-\eta_{21}) \wedge (4\,\eta_4-\eta_5)
-(\eta_7-\eta_9) \wedge (3\,\sigma_{13}+2\,\sigma_{22}+2\,\sigma_{33})
+\eta_8 \wedge (2\, \sigma_{23}-3\,\sigma_{34})
\nonumber \\
\fl
&&
+(\eta_{10}-\eta_{15}) \wedge (2\,\sigma_{23}-3\,\sigma_{34})
-(\eta_{13}-\eta_{16}) \wedge (3\,\sigma_{13}+2\,\sigma_{22}+3\,\sigma_{33})
\nonumber \\
\fl
&&
-\eta_{19} \wedge (3\,\sigma_{13}+2\,\sigma_{22})
+\eta_{20} \wedge (4\,\eta_3+2\,\eta_6-5\,\sigma_{23}),
\nonumber 
\\
%\end{eqnarray}
%%%%%%%%%%%%%%%%%%%%%%%%%%%%%%%%%%%%%%%%%%%%%%%%%%%%%%%%%%%%%%%%%%%%%%%%%%%%%%%%%%%%
%\begin{eqnarray}
%------------------------------------------------------------------------
\fl
d\eta_{22}  &=&  
(\eta_{24}+\eta_{25}-\eta_{36}+\eta_{33}+\eta_{39}) \wedge \theta_1
+(2\,(\eta_{23}-\eta_{26}+\eta_{35})+\eta_{40}) \wedge \theta_2
-3\,(\eta_7 +\eta_{19}) \wedge \sigma_{14}
\nonumber \\
\fl
&&
+(\eta_{25}+\eta_{33}-\eta_{36}) \wedge \theta_3
-(\eta_{23}-\eta_{26}+\eta_{35}) \wedge \theta_4
+\eta_{38} \wedge \xi^1
+(3\,\eta_{18}+2\,\eta_{21}) \wedge \eta_2
\nonumber \\
\fl
&&
+(2\,\eta_{23}-\eta_{26}+2\,\eta_{35}+\eta_{37}+\eta_{40}) \wedge \xi^2
+\eta_{41} \wedge \xi^3
+\eta_{42} \wedge \xi^4
+3\,\eta_{21} \wedge (2\,\eta_3+\eta_6-2\,\sigma_{23})
\nonumber \\
\fl
&&
-2\,(\eta_8 + 2 \,(\eta_{10}-\eta_{15})-\eta_{17}-\eta_{20}) \wedge \eta_1
+3\,(\eta_9-\eta_{13}+\eta_{16}) \wedge (\sigma_{14}+\sigma_{34})
\nonumber \\
\fl
&&
+(4\,(\eta_8+2\,(\eta_{10}-\eta_{15})-\eta_{17}-\eta_{20})+3\,\eta_{22}) \wedge \eta_4
+3\,\eta_{17} \wedge \sigma_{33}
-3\,\eta_{20} \wedge (\sigma_{13}+\sigma_{22})
\nonumber \\
\fl
&&
-(\eta_8+2\,(\eta_{10}-\eta_{15})-\eta_{17}-\eta_{20}+\eta_{22}) \wedge \eta_5
+3\,\eta_8 \wedge (2\,(\sigma_{13}+\sigma_{22})+\sigma_{33}-\sigma_{44})
\nonumber \\
\fl
&&
+3\,(\eta_{10}-\eta_{15})\wedge (2\, \sigma_{13}+2\,\sigma_{22}-\sigma_{44}).
\nonumber
\end{eqnarray}

\end{document}